\documentclass{amsart}
\usepackage{amsfonts}
\usepackage{bbm}
\usepackage{amsfonts,amssymb,amsmath,amsthm}
\usepackage{url}
\usepackage{enumerate}
\usepackage{bbm}
\usepackage[all]{xy}

\newtheorem{thrm}{Theorem}[section]
\newtheorem{lem}[thrm]{Lemma}

\newtheorem{cor}[thrm]{Corollary}
\theoremstyle{definition}
\newtheorem{definition}[thrm]{Definition}

\newtheorem{nota}[thrm]{Notation}
\numberwithin{equation}{section}

\author{Xia Zhao}
\address{School of Mathematical Sciences, Tongji University, Shanghai, China \ 200092}
\email{1710383@tongji.edu.cn}

\author{Xiaochun Fang}
\address{School of Mathematical Sciences, Tongji University, Shanghai, China \ 200092}
\email{xfang@tongji.edu.cn}

\author{Qingzhai Fan}
\address{Department of Mathematics, Shanghai Maritime University, Shanghai, China \ 201306}
\email{qzfan@shmtu.edu.cn}

\begin{document}

\title[Some Permanence for Large Subalgebra ]{Some permanence for Large subalgebra}

\begin{abstract}
In this paper, we give two properties of C*-algebra that could be deduced from the properties of its large subalgebra. Let $A$ be an infinite dimensional simple unital C*-algebra and let $B$ be a centrally large subalgebra of $A$, we prove that $A$ has real rank zero if $B$ has real rank zero. If $A$ is stablely finite in addition, $B$ is a large subalgebra of $A$, we prove that $B$ has local weak comparison if $A$ has local weak comparison, and $A$ has local weak comparison if $M_{2}(B)$ has local weak comparison. As a consequence, we show that $A$ has weak comparison if and only if $B$ has weak comparison. These results could be used to study some properties of C*-algebra from its large subalgebra or centrally large subalgebra.

\end{abstract}

\maketitle
\section{introduction}
Putnam \cite{I. F. P1} introduced Putnam subalgebra of the crossed product by a minimal homeomorphism. Let $X$ be the Cantor set, and let $h:X\rightarrow X$ be a minimal homeomorphism. Let $u$ be the standard unitary in the crossed product $C^{*}(\mathbb{Z},X,h)$, and fix $y\in X$. Then the Putnam subalgebra of $C^{*}(\mathbb{Z},X,h)$ is generated by $C(X)$ and all elements $fu$ with $f\in C(X)$ satisfying $f(y)=0$. It is closely enough to use information about the Putnam subalgebra to obtain information about $C^{*}(\mathbb{Z},X,h)$. There are so many applications about this in papers \cite{Q. L1}, \cite{G. A. E1}, \cite{N. C. P1} etc. .

Phillips \cite{N. C. P2} and Archey etc. \cite{N. C. P3} gave the definitions of large subalgebra and centrally large subalgebra which are abstractions of Putnam subalgebra.  An basic example of large subalgebra was presented in \cite{N. C. P2}. Let $X$ be an infinite compact metric space and let $h:X\rightarrow X$ be a minimal homeomorphism. Suppose $Y\subset X$ is closed. The $Y$-$orbit \, breaking \, subalgebra$ of $C^{*}(\mathbb{Z},X,h)$ associated to $Y$ is the subalgebra generated by $C(X)$ and all elements $fu$ with $f\in C(X)$ satisfying $f|_{Y}=0$. If $Y$ meets each orbit at most once, then this algebra is (centrally) large subalgebra in $C^{*}(\mathbb{Z},X,h)$.

 Large subalgebra is much easier than original algebra. So to study the properties of C*-algebra $A$, we usually study the properties of large subalgebra $B$ at first. Now there is a question that which properties of $B$ could be used to deduce properties of $A$. Phillips \cite{N. C. P2} proved the following results: $B$ and $A$ have the same radius of  comparison when $A$ is stably finite; $B$ is finite if and only if $A$ is finite; $B$ is pure infinite if and only if $A$ is pure infinite; the restriction maps $T(A)\rightarrow T(B)$ and $QT(A)\rightarrow QT(B)$ (on tracial states and quasitraces) are bijective. Archey etc. \cite{N. C. P3} showed $A$ is stable rank one if  $B$ is stable rank one, $A$ is real rank zero if $B$ is real rank zero and stable rank one. Archey etc. \cite{N. C. P4} got the permanence of tracially $\mathcal{Z}$-absorption and $\mathcal{Z}$-absorption in some conditions.

In addition to the above properties, what other property of a C*-algebra can be deduced from the property of its large subalgebra? Phillips \cite{N. C. P5} gave some open problems about this question, like tracial rank zero, finite nuclear dimension etc. . Question 1.47 in \cite{N. C. P5} asked that if $B$ is real rank zero, does it follow that $A$ is real rank zero without the condition that $B$ is stable rank one. And Archey etc. \cite{N. C. P3} said it is possible that one does not need the subalgebra to have stable rank one. Our paper gives an affirmative answer for this problem.

Besides, comparison is an important property of C*-algebra. Toms and Winter \cite{A. T1} conjecture that strict comparison of positive elements, finite nuclear dimension and $\mathcal{Z}$-stable are equivalent in infinite dimensional unital nuclear separable C*-algebra.  \cite{N. C. P5} proved $B$ and $A$ have the same comparison radius. It is called strict comparison when the comparison radius is zero, that is, $B$ has strict comparison if and only if $A$ has strict comparison. However, there are many other comparison properties. Ortega etc. \cite{E. O1} gave a definition of $n$-comparison of Cuntz semigroup and Winter \cite{W. W1} introduced $m$-comparison of separable simple unital C*-algebra. The above two definitions are in fact equivalent for simple C*-algebra. Kirchberg and R{\o}rdam defined local weak comparison and weak comparison in \cite{E. K1}. They proved local weak comparison and weak comparison are weaker than $m$-comparison and strict comparison in some C*-algebra. Fan etc. \cite{Q. F1} proved the inheritance of $m$-comparison from its large subalgebra. In this paper, we prove the permanence of local weak comparison and weak comparison.

To be precise, we first study that real rank zero of centrally large subalgebra could deduce real rank zero of the original algebra without the condition of stable rank one. We get the following result.

(1) Let $A$ be an infinite dimensional simple unital C*-algebra, and let $B\subset A$ be a centrally large subalgebra. If $RR(B)=0$, then we have $RR(A)=0$.

  With the property of stable rank one, Archey etc. \cite{N. C. P3} used matrix decomposition by positive elements to prove the permanence of stable rank one. It is not nice about the Cuntz comparison and Cuntz semigroup of $B$ when $B$ is only  real rank zero, we can not use the same proof in \cite{N. C. P3}. However, C*-algebra with real rank zero has property (SP), there are so many projections in $B$. So we use matrix decomposition by projections to prove this result.

Next for the permanence of local weak comparison and weak comparison, we have the following results.

Let $A$ be an infinite dimensional stably finite simple separable unital C*-algebra, and let $B\subset A$ be a large subalgebra.

 (1) $B$ has local weak comparison if $A $ has  local weak comparison.

 (2) $A$ has  local weak comparison if $M_{2}(B)$ has local weak comparison.

 (3) $A$ has  weak comparison if and only if $B$ has weak comparison.

  Since local weak comparison of $A$ and $M_{2}(A)$ are not the same, we could prove $B$ has  local weak comparison when $A $ has  local weak comparison. But for the other direction, we only get $A$ has  local weak comparison if $M_{2}(B)$ has local weak comparison. (3) is a nature corollary from (1) and (2).

The paper is organized as follows. Section 2 contains some preliminaries about Cuntz subequivalent and large subalgebra. Section 3
presents real rank zero of centrally large algebra could deduce real rank zero of the original C*-algebra. Section 4 shows the permanence of local weak comparison and weak comparison.

\section{preliminaries}
In this  section, we introduce some definitions, symbols and known facts about Cuntz subequivalent and large subalgebra.

For a C*-algebra $A$, let $M_{\infty}(A)$ denote the algebraic direct limit of system $(M_{n}(A))_{n=1}^{\infty}$,
$K\otimes A$ denote the C*-algebraic direct limit of system $(M_{n}(A))_{n=1}^{\infty}$ and $A_{+}$ denote the set of all positive elements in $A$.

\begin{definition}
Let $A$ be a C*-algebra, $a\in A_{+}$ and $\varepsilon>0$. Let $f:[0,\infty)\rightarrow[0,\infty)$ be the function

$$f(\lambda)=(\lambda-\varepsilon)_{+}=\begin{cases}
              0 & 0\leq\lambda\leq\varepsilon \\
              \lambda-\varepsilon &\varepsilon<\lambda.\\
             \end{cases}$$
Then define $(a-\varepsilon)_{+}=f(a).$
\end{definition}

\begin{definition}
Let $A$ be a C*-algebra and $a,b\in(K \otimes A)_{+}$.

(1) We say that $a$ is $Cuntz\,subequivalent$ to $b$ over $A$, written $a\precsim_{A}b$, if there is a sequence $(v_{n})_{n=1}^{\infty}$ in $K\otimes A$ such that $lim_{n\rightarrow\infty}v_{n}bv_{n}^{*}=a$.

(2) We say that $a$ and $b$ are $Cuntz\,equivalent$ in $A$, written $a\sim_{A}b$, if $a\precsim_{A}b$ and $b\precsim_{A}a$. We write $\langle a\rangle$ for the equivalence class of $a$.

(3) The $Cuntz\,semigroup$ of $A$ is
\[{\rm Cu}(A)=(K\otimes A)_{+}/\sim_{A},\]
together with the operation $\langle a\rangle+\langle b\rangle=\langle a\oplus b\rangle$ and the partial order $\langle a \rangle\leq\langle b \rangle\Leftrightarrow a\precsim_{A}b$.

(4) We define the semigroup
\[W(A)=M_{\infty}(A)_{+}/\sim_{A}\]
with the same operation.
\end{definition}

Some known facts about Cuntz subequivalent are in the following lemma. All the proofs could be found in section 1 of \cite{N. C. P2}.

\begin{lem}\label{fundamental}
Let $A$ be a C*-algebra.

{\rm(1)} Let $a,\,b\in A_{+}$ satisfy $a\leq b$, then $a\precsim_{A}b$.

{\rm(2)} Let $a\in A_{+}$ and $\varepsilon_{1},\,\varepsilon_{2}>0$, then $((a-\varepsilon_{1})_{+}-\varepsilon_{2})_{+}=(a-(\varepsilon_{1}+\varepsilon_{2}))_{+}$.

{\rm(3)} Let $\varepsilon>0$, $a,\,b\in A_{+}$ satisfy $\|a-b\|<\varepsilon$, then $(a-\varepsilon)_{+}\precsim_{A}b$.

{\rm(4)} Let $a,b\in A_{+}$, then the following are equivalent:

$\quad\quad${\rm(a)} $a\precsim_{A}b$;

$\quad\quad${\rm(b)} $(a-\varepsilon)_{+}\precsim_{A}b$ for all $\varepsilon>0$;

$\quad\quad${\rm(c)} For every $\varepsilon>0$, there is $\delta>0$ such that $(a-\varepsilon)_{+}\precsim_{A}(b-\delta)_{+}$.

{\rm(5)} Let $\varepsilon>0$, $\lambda>0$ and $a,\,b\in A$ satisfy $\|a-b\|<\varepsilon$, then  $(a-\lambda-\varepsilon)_{+}\precsim_{A}(b-\lambda)_{+}$.

{\rm(6)} Let $a,\,g\in A_{+}$ with $0\leq g\leq1$, and let $\varepsilon>0$. Then
$$(a-\varepsilon)_{+}\precsim_{A}[(1-g)a(1-g)-\varepsilon]_{+}\oplus g.$$
\end{lem}

\begin{nota}\label{quasitraces} Let $A$ be a unital C*-algebra.

(1) We denote $QT(A)$ the set of normalized 2-quasitraces on $A$ (Definition 2.31 of \cite{P. A1} and definition 2.1.1 of \cite{B. B1}).

(2) Define $d_{\tau}:M_{\infty}(A)_{+}\rightarrow[0,\infty)$ by $d_{\tau}(a)=\lim_{n\rightarrow\infty}\tau(a^{\frac{1}{n}})$ for all $a\in M_{\infty}(A)_{+}$ and $\tau\in QT(A)$. We also use the same notation for the corresponding functions on $(K\otimes A)_{+}$, Cu$(A)$ and $W(A)$. $d_{\tau}$ is well defined on Cu$(A)$ and $W(A)$ by part of the proof of proposition 4.2 in \cite{G. A. E2}. It follows that $d_{\tau}$ defines a state on $W(A)$ by the proof of theorem 2.32 in \cite{P. A1}
\end{nota}

The following result is well known and its proof could be found in lemma 1.22 in \cite{N. C. P2}.
\begin{lem}\label{lower}\cite[Lemma 1.22]{N. C. P2}
Let $A$ be a unital C*-algebra, and let $a\in(K\otimes A)_{+}$. Then the function $\tau\rightarrow d_{\tau}(a)$ is lower semicontinuous from $QT(A)$ to $[0,\infty]$.
\end{lem}

First, we recalled  the (centrally) large subalgebra defined by Phillips in \cite{N. C. P2} and Archey etc. \cite{N. C. P3}.

\begin{definition}\label{large subalgebra}
Let $A$ be an infinite dimensional simple unital C*-algebra. A unital subalgebra $B\subset A$ is said to be $large$ in $A$ if for every $m\in \mathbb{Z}_{+},\,a_{1},\,a_{2},\,...,\,a_{m}\in A,\,\varepsilon>0,\,x\in A_{+}$ with $\| x\|=1$, and $y\in B_{+}\setminus\{0\}$, there are $c_{1},\,c_{2},\,...,\,c_{m}\in A$ and $g\in B$ such that:

(1) $0\leq g\leq1;$

(2) For $j=1,\,2,\,...,\,m$, we have $\|c_{j}-a_{j}\|<\varepsilon;$

(3) For $j=1,\,2,\,...,\,m$, we have $(1-g)c_{j}\in B;$

(4) $g\precsim_{B}y,\,g\precsim_{A}x;$

(5) $\|(1-g)x(1-g)\|>1-\varepsilon;$

\noindent We say that $B$ is $centrally\,\, large$ in $A$ if we require that in addition:

(6) For $j=1,\,2,\,...,\,m$, we have $\|ga_{j}-a_{j}g\|<\varepsilon.$
\end{definition}

In definition \ref{large subalgebra}, the elements $c_{1},\,c_{2},\,...,\,c_{m}$ could be chosen such that $\|c_{j}\|\leq\|a_{j}\|$ for $j=1,\,2,\,...,\,m$. If the elements $a_{1},\,a_{2},\,...,\,a_{m}\in A_{+}$, then $c_{1},\,c_{2},\,...,\,c_{m}$ could be chosen in $A_{+}$. The following lemma is from lemma 4.7 of \cite{N. C. P2} and lemma 3.4 of \cite{N. C. P3}.

\begin{lem}\label{def}
 Let $A$ be an infinite dimensional simple unital C*-algebra and $B\subset A$ be a large subalgebra. Let $m,\,n\in \mathbb{Z}_{+}$, $a_{1},\,a_{2},\,...,\,a_{m}\in A,\,b_{1},\,b_{2},\,...,\,b_{n}\in A_{+},\,\varepsilon>0,\,x\in A_{+}$ with $\| x\|=1$, and $y\in B_{+}\setminus\{0\}$. Then there are $c_{1},\,c_{2},\,...,\,c_{m}\in A,\,d_{1},\,d_{2},\,...,\,d_{n}\in A_{+}$ and $g\in B$ such that:

{\rm(1)} $0\leq g\leq1;$

{\rm(2)} $\|c_{j}-a_{j}\|<\varepsilon$ for $j=1,\,2,\,...,\,m$, $\|b_{i}-d_{i}\|<\varepsilon$ for $i=1,\,2,\,...,\,n$;

{\rm(3)} $\|c_{j}\|\leq \| a_{j}\|$ for $j=1,\,2,\,...,\,m$, $\| b_{i}\|\leq \| d_{i}\|$ for $i=1,\,2,\,...,\,n$;

{\rm(4)} $(1-g)c_{j}\in B$ for $j=1,\,2,\,...,\,m$, $(1-g)d_{i}(1-g)\in B$ for $i=1,\,2,\,...,\,n$;

{\rm(5)} $g\precsim_{B}y,\,g\precsim_{A}x;$

{\rm(6)} $\| (1-g)x(1-g)\|>1-\varepsilon$;

\noindent If $B\subset A$  is a centrally large subalgebra, then it could require in addition:

{\rm(7)} $\|ga_{j}-a_{j}g\|<\varepsilon$ for $ j=1,\,2,\,...,\,m$, $\|gb_{i}-b_{i}g\|<\varepsilon$ for $i=1,\,2,\,...,\,n$.
\end{lem}

\section{Real Rank Zero}

In this section,  we prove that $A$ is real rank zero if $B$ is real rank zero when $B\subset A$ is a centrally large subalgebra. First we introduce some properties of C*-algebra with real rank zero.
\begin{definition}\label{real rank zero}
A unital C*-algebra is said to have $real\,\, rank\,\, zero$, written $RR(A)=0$, if the set of invertible self-adjoint elements is dense in $A_{sa}$. A non-unital C*-algebra is said to have $real\,\, rank\,\, zero$ if $RR(\widetilde{A})=0$.
\end{definition}
\begin{definition}\label{sp}
We say a C*-algebra $A$ has $property\,\, (SP)$ if every nonzero hereditary C*-algebra of $A$ contains a nonzero projection.
\end{definition}
It is clear that every C*-algebra with real rank zero has property (SP).
\begin{lem}\label{property}\cite[Lemma 3.5.6]{H. X. L1}
Let $A$ be a simple C*-algebra with property (SP) and $p\in A$ be a non zero projection. Suppose that $a\in A_{+}$ is a nonzero element. Then there is a nonzero projection $q\in\overline{aAa}$ such that $q\precsim p$.
\end{lem}
\begin{lem}\label{decomposition}\cite[Lemma 3.5.7]{H. X. L1}
Let $A$ be a non-elementary simple C*-algebra with property (SP). Then for any nonzero projection $p\in A$ and any integer $n>0$, there are $n+1$ mutually orthogonal projections $q_{1},\,q_{2},\,...,\,q_{n+1}$ such that $q_{1}\neq0,\,q_{1}\sim q_{i},\,i=1,\,2,\,...,\,n$, and $p=q_{1}+q_{2}+...+q_{n+1}$.
\end{lem}
Next, we give some lemmas needed in the proof of our theorem and some proofs of the lemmas could be found in \cite{N. C. P2} and \cite{N. C. P3}.
\begin{lem}\label{squares}\cite[Lemma 2.5]{N. C. P3}
 Suppose that $f:[0,1]\rightarrow \mathbb{C}$ is continuous. Then for every $\varepsilon>0$, there is $\delta>0$ such that whenever $A$ is a C*-algebra and $x,\,y\in A$ satisfying
 $$0\leq x\leq1,\qquad \| y\|\leq1\qquad and \qquad \| xy-yx\|<\delta,$$
 then $\| f(x)y-yf(x)\|<\varepsilon$.
 \end{lem}
The following lemma is slightly different from lemma 2.6 of \cite{N. C. P3}.
\begin{lem}\label{dist} Suppose $f:[0,1]\rightarrow [0,1]$ is continuous with $f(0)=0$. Then for every $\varepsilon>0$, there is $\delta>0$ satisfying the following.

  Let $A$ be a C*-algebra, $B\subset A$ be a subalgebra, and let $x\in A_{+}$, $a\in B$ with $\|x\|\leq1,\,0\leq a\leq1$, there is $b\in B_{+}$ such that $\| axa-b\|<\delta$. Then there exists $c\in B_{+}$ such that $\|f(a)xf(a)-c\|<\varepsilon$.
 \end{lem}
\begin{proof}

 Given $\varepsilon>0$, we may assume $1>\varepsilon>0$, then there exist $n\in \mathbb{Z}_{\geq0}$ and a polynomial $g(\lambda)=\sum_{k=0}^{n}\alpha_{k}\lambda^{k}$ with $\alpha_{k}\in\mathbb{R}$ for $k=0,1,...,n$ such that $\alpha_{0}=0$ and $\mid g(\lambda)-f(\lambda)\mid<\frac{\varepsilon}{4}$.

 Define $\delta=\frac{\varepsilon}{4\sum_{j,k=1}^{n}|\alpha_{j}||\alpha_{k}|+1}$. Let $x\in A$ and $a\in B$ be as in the hypotheses. Choose $b\in B_{+}$ such that $\|axa-b\|<\delta$. Set $c=(\sum_{k=1}^{n}\alpha_{k}a^{k-1})b(\sum_{k=1}^{n}\alpha_{k}a^{k-1})$, then
  \begin{align}
&\|g(a)xg(a)-c\|\notag\\
&=\| (\sum_{k=1}^{n}\alpha_{k}a^{k})x(\sum_{k=1}^{n}\alpha_{k}a^{k})-(\sum_{k=1}^{n}\alpha_{k}a^{k-1})b(\sum_{k=1}^{n}\alpha_{k}a^{k-1})\|\notag\\
&=\| \sum_{k,j=1}^{n}\alpha_{k}\alpha_{j}a^{k-1}(axa-b)a^{j-1}\|\notag\\
&\leq \sum_{k,j=1}^{n}|\alpha_{k}||\alpha_{j}|\|axa-b\|\|a^{k-1}\|\|a^{j-1}\|\notag\\
&\leq \sum_{k,j=1}^{n}|\alpha_{k}||\alpha_{j}|\delta<\frac{\varepsilon}{4}\notag.
\end{align}
Since $\| f(a)\|\leq1,\| g(a)\|\leq1+\frac{\varepsilon}{4}$, we have
 \begin{align}
&\| f(a)xf(a)-c\|\notag\\
&\leq\| f(a)\|\| x\|\| f(a)-g(a)\|+\| f(a)-g(a)\|\| x\|\|g(a)\|+\| g(a)xg(a)-c\|\notag\\
&\leq\frac{\varepsilon}{4}+(1+\frac{\varepsilon}{4})\frac{\varepsilon}{4}+\frac{\varepsilon}{4}\notag\\
&<\varepsilon\notag.
\end{align}
\end{proof}

\begin{lem}\label{large}\cite[Lemma 5.3]{N. C. P2}
 Let $A$ be an infinite dimensional simple unital C*-algebra, and let $B\subset A$ be a large subalgebra. Let $r\in B_{+}\setminus\{0\}$ and $a\in \overline{rAr}$ be positive with $\|a\|=1$. For $\varepsilon>0$, there is a positive element $b\in\overline{rBr}$ such that:

 (1) $\| b\|=1 $;

 (2) $b\precsim_{A}a$;

 (3) $\|ab-b\|<\varepsilon $.

\end{lem}

\begin{lem}\label{continuous}\cite[Lemma 2.5.11]{H. X. L1}
 Let $f\in C([0,1])$ with $f(0)=0$. For every $\varepsilon>0$, there exists $\delta>0$ satisfying the following: suppose that $A$ is a C*-algebra, $a,\,b\in A_{+}$ with $\|a\|\leq1,\,\| b\|\leq1$. If $\|ab-b\|<\delta$, then $\|af(b)-f(b)\|<\varepsilon$.
 \end{lem}

\begin{lem}\label{approximate}
Let $A$ be an infinite dimensional simple unital C*-algebra, and let $B\subset A$ be a large subalgebra. For any $\varepsilon>0$, any $x\in A_{sa}$ with $\| x\|\leq1$, $a\in A_{+}\setminus\{0\}$ with $ax=xa=0$. Then there exist $y\in A_{sa}$ with $\| y\|\leq1$ and $b\in B_{+}\setminus\{0\}$  such that $\|x-y\|<\varepsilon$ and $by=yb=0$.
\end{lem}
\begin{proof}
We may assume $\|a\|=1$. Define $f_{0},\,f_{1}:[0,1]\rightarrow [0,1]$ by

 $$f_{0}(\lambda)=\begin{cases}
                 0 & 0\leq\lambda\leq\frac{1}{4} \\
              4x & \frac{1}{4}\leq\lambda\leq \frac{1}{2} \\
              1 & \frac{1}{2} \leq\lambda\leq1,\\
                 \end{cases}$$
$$f_{1}(\lambda)=\begin{cases}
                 0 & 0\leq\lambda\leq\frac{1}{2} \\
              4x-2 & \frac{1}{2}\leq\lambda\leq \frac{3}{4} \\
              1 & \frac{3}{4}\leq\lambda\leq 1.\\
                 \end{cases}$$

 Then $f_{0}f_{1}=f_{1}$. For $f_{0}\in C([0,1])$ and $\frac{\varepsilon}{2}$, by lemma \ref{continuous}, there is $\delta>0$ such that $\|af_{0}(b)-f_{0}(b)\|<\frac{\varepsilon}{2}$ if $\|ab-b\|<\delta$.

    For $\delta>0$, $a\in A_{+}$ with $\|a\|=1$, since $B\subset A$ is a large subalgebra, by lemma \ref{large}, we  have $b\in B_{+}$ such that
$\|b\|=1, \|ab-b\|<\delta $. Then $\|af_{0}(b)-f_{0}(b)\|<\frac{\varepsilon}{2}$.

  Let $r_{0}=f_{0}(b),r_{1}=f_{1}(b)$, then we have

  (1) $0\leq r_{0}\leq1,0\leq r_{1}\leq1,\|r_{0}\|=1,\|r_{1}\|=1$;

 (2) $r_{0}r_{1}=r_{1}$;

 (3) $\|a r_{0}-r_{0}\|<\frac{\varepsilon}{2}$.

  Let $y=(1- r_{0})x(1- r_{0}),\,b=r_{1}$. Then $y\in A_{sa}$ with $\|y\|\leq1$, $b\in B_{+}$, and $by=r_{1}y=0,\,yb=yr_{1}=0$.

  Since $\|x r_{0}\|\leq\| x r_{0}-x ar_{0}\|\leq\| r_{0}-ar_{0}\|<\frac{\varepsilon}{2}$,
  similarly, $\| r_{0}x\|<\frac{\varepsilon}{2}$. Then
  \begin{align}
  \|x-y\|&=\|xr_{0}+r_{0}x(1-r_{0})\|\notag\\
   & \leq\| xr_{0}\|+\| r_{0}x(1-r_{0})\|\notag\\
   & <\frac{\varepsilon}{2}+\frac{\varepsilon}{2}=\varepsilon\notag.
\end{align}
\end{proof}

\begin{lem}\label{projection}
 Let $A$ be an infinite dimensional simple unital C*-algebra, and let $B\subset A$ be a centrally large subalgebra. $p$ is any projection in $B$, supose $m,\,n\in \mathbb{Z}_{\geq0},\,\{pa_{1}p,\,pa_{2}p,\,...,\,pa_{m}p\}\subset pAp$, $\{pb_{1}p,\,pb_{2}p,\,...,\,pb_{n}p\}\subset (pAp)_{+},\,\varepsilon>0,\,x\in (pAp)_{+} $ with $\| x \|=1$, $y\in (pBp)_{+}\setminus\{0\}$, then there are $\{c_{1},\,c_{2},\,...,\,c_{m}\}\subset pBp$, $\{d_{1},\,d_{2},\,...,\,d_{n}\}\subset (pBp)_{+}$, $g\in pBp$ such that:

 {\rm(1)} $0\leq g\leq 1$;

 {\rm(2)} $\|(p-g)pa_{j}p-c_{j}\|<\varepsilon$ for $j=1,\,2,\,...,\,m,$ $\|(p-g)pb_{i}p(p-g)-d_{i}\|<\varepsilon$ for $i=1,\,2,\,...,\,n$;

 {\rm(3)} $g\precsim_{pBp}y,\, g\precsim_{pAp}x$;

 {\rm(4)} $\|(p-g)x(p-g)\|>1-\varepsilon$;

 {\rm(5)} $\| gpa_{j}p-pa_{j}pg\|<\varepsilon$  for $j=1,\,2,\,...,\,m,$ $\|gpb_{i}p-pb_{i}pg\|<\varepsilon$ for $i=1,\,2,\,...,\,n$.

\end{lem}
\begin{proof}
Let $F=\{pa_{1}p,\,pa_{2}p,\,...,\,pa_{m}p,\,p\}\subset A$, $G=\{pb_{1}p,\,pb_{2}p,\,...,\,pb_{n}p\}\subset A_{+}$, $x\in A_{+} $ with $\| x \|=1$ and $y\in B_{+}\setminus\{0\}$. Since $B$ is centrally large subalgebra of $A$, by lemma \ref{def}, there exist $e_{1},\,e_{2},\,...,\,e_{m},\,d\in A$, $f_{1},\,f_{2},\,...,\,f_{n} \in A_{+}$ and $g_{1}\in B$ such that:

 $(1)^{'}$ $0\leq g_{1}\leq 1$;

 $(2)^{'}$ $\| e_{j}-pa_{j}p\|<\frac{\varepsilon}{3},\,\|d-p\|<\frac{\varepsilon}{3}$ for $j=1,\,2,\,...,\,m,$ $\|f_{i}-pb_{i}p\|<\frac{\varepsilon}{3}$ for $i=1,\,2,\,...,\,n$;

 $(3)^{'}$ $\|e_{j}\|\leq \|pa_{j}p\|$ for $j=1,\,2,\,...,\,m$, $\| f_{i}\|\leq \|pb_{j}p\|$ for $i=1,\,2,\,...,\,n$;

$(4)^{'}$ $(1-g_{1})e_{j}\in B$ for $j=1,\,2,\,...,\,m$, $(1-g_{1})f_{i}(1-g_{1})\in B$ for $i=1,\,2,\,...,\,n$;

 $(5)^{'}$ $g_{1}\precsim_{B}y, \,g_{1}\precsim_{A}x$;

 $(6)^{'}$ $\|(1-g_{1})x(1-g_{1})\|>1-\frac{\varepsilon}{3}$;

 $(7)^{'}$ $\| g_{1}pa_{j}p-pa_{j}pg_{1}\|<\frac{\varepsilon}{3}$ for $j=1,\,2,\,...,\,m$, $\| g_{1}p-pg_{1}\|<\frac{\varepsilon}{3}$, $\|g_{1}pb_{i}p-pb_{i}pg_{1}\|<\frac{\varepsilon}{3}$ for $i=1,2,...,n$.

 Let $g=pg_{1}p$, then (1) $0\leq g\leq 1$ holds.

 (2). Let $c_{j}=p(1-g_{1})e_{j}p$ for $j=1,2,...,m$, $d_{i}=p(1-g_{1})f_{i}(1-g_{1})p$ for $i=1,2,...,n$, then $c_{j}\in pBp$, $d_{i}\in (pBp)_{+}$ and

 $\|(p-g)pa_{j}p-c_{j}\|=\|(p-g)pa_{j}p-p(1-g_{1})e_{j}p\|<\varepsilon$,

  $\|(p-g)pb_{i}p(p-g)-d_{i}\|=\|(p-g)pb_{i}p(p-g)-p(1-g_{1})f_{i}(1-g_{1})p\|<\varepsilon$.

 So there exist $c_{j},\,d_{i}$ such that $\|(p-g)pa_{j}p-c_{j}\|<\varepsilon$ for $j=1,\,2,\,...,\,m,$ $\|(p-g)pb_{i}p(p-g)-d_{i}\|<\varepsilon$ for $i=1,\,2,\,...,\,n$.

 (3). Since $g_{1}\precsim_{B}y,\, g_{1}\precsim_{A}x$, then there exist $\{r_{k}\}\subset A,\,\{c_{k}\}\subset B$ such that $r_{k}xr_{k}^{*}\rightarrow g_{1}, \,c_{k}yc_{k}^{*}\rightarrow g_{1}$. Then $\{pr_{k}p\}\subset pAp,\,\{pc_{k}p\}\subset pBp$, and $pr_{k}xr_{k}^{*}p\rightarrow pg_{1}p, pc_{k}yc_{k}^{*}p\rightarrow pg_{1}p$. Thus $g\precsim_{pBp}y,\, g\precsim_{pAp}x$.

(4). Since $\| g_{1}p-pg_{1}\|<\frac{\varepsilon}{3}$, then
\begin{align}
  &\|(p-g)x(p-g)\|\notag\\
  &= \|p(1-g_{1})pxp(1-g_{1})p\|\notag\\
   &= \| p(1-g_{1})pxp(1-g_{1})p-(1-g_{1})x(1-g_{1})+(1-g_{1})x(1-g_{1})\|\notag\\
   &\geq\|(1-g_{1})x(1-g_{1})\|-\| p(1-g_{1})pxp(1-g_{1})p-(1-g_{1})x(1-g_{1})\|\notag\\
   & >\|(1-g_{1})x(1-g_{1})\|-\frac{2\varepsilon}{3}>1-\varepsilon\notag.
\end{align}

(5). Since

$\|gpa_{j}p-pa_{j}pg\|=\| pg_{1}ppa_{j}p-pa_{j}ppg_{1}p\|\leq\| g_{1}pa_{j}p-pa_{j}pg_{1}\|<\varepsilon,$

$ \| gpb_{i}p-pb_{i}pg\|=\|pg_{1}ppb_{i}p-pb_{i}ppg_{1}p\|\leq\|g_{1}pb_{i}p-pb_{i}pg_{1}\|<\varepsilon.$

 Then $\|gpa_{j}p-pa_{j}pg\|<\varepsilon$ for $j=1,\,2,\,...,\,m$, $ \| gpb_{i}p-pb_{i}pg\|<\varepsilon$ for $i=1,\,2,\,...,\,n$.

\end{proof}

\begin{thrm} Let $A$ be an infinite dimensional simple unital C*-algebra, and let $B\subset A$ be a centrally large subalgebra. If $RR(B)=0$, we have $RR(A)=0$.
\end{thrm}

\begin{proof} Let $x$ be a self-adjoint element in $A$ with $\|x\|=1$, and let $1>\varepsilon>0$. We will show that there is an invertible self-adjoint element $z\in A$ such that $\| x-z\|<\varepsilon$. We decompose our proof in the following three steps.

\vspace{10pt}
Step $\uppercase\expandafter{\romannumeral1}$, we prove that there are projection $p\in B_{+}$ and self-adjoint element $y\in (1-p)A(1-p)$ with $\|y\|\leq1$ such that $\| y-x\|<\frac{\varepsilon}{3}$.

By continuous function calculus, first we get  $x_{0}\in A_{sa}$ with $\| x_{0}\|\leq1$ and $y_{0}\in A_{+}\setminus\{0\}$ such that $\| x_{0}-x\|<\frac{\varepsilon}{6}$ and $x_{0}y_{0}=y_{0}x_{0}=0$.

Let $f(t)\in C([-1,1])$ with $0\leq f\leq1$, $f(t)=1$ if $\mid t\mid<\frac{\varepsilon}{24}$ and $f(t)=0$ if $\mid t\mid\geq\frac{\varepsilon}{12}$.
If $0\notin sp(x)$, then $x$ is invertible. So we assume that $0\in sp(x)$, thus $f(x)\neq0$.
Let $f_{0}(t)\in C([-1,1])$ with $0\leq f_{0}\leq1$, $f_{0}(t)=1$ if $\mid t\mid>\frac{\varepsilon}{6}$ and $f_{0}(t)=0$ if $\mid t\mid\leq\frac{\varepsilon}{12}$. Let $x_{0}=f_{0}(x)x,\,y_{0}=f_{0}(x)$, then $x_{0}\in A_{sa}$ with $\| x_{0}\|\leq1$ and $y_{0}\in A_{+}\setminus\{0\}$ such that $x_{0}y_{0}=y_{0}x_{0}=0$ and $\| x_{0}-x\|<\frac{\varepsilon}{6}$.

 By lemma \ref{approximate}, there exist $y\in A_{sa}$ with $\| y\|\leq1$ and $b\in B_{+}$ such that $\| y-x_{0}\|<\frac{\varepsilon}{6}$ and $by=yb=0$. Then $\| y-x\|\leq\| y-x_{0}\|+\| x_{0}-x\|<\frac{\varepsilon}{3}$. Since $RR(B)=0$, then $B$ has property (SP), so there is a nonzero projection $p\in \overline{bBb}$ such that $py=yp=0$, then $y\in (1-p)A(1-p)$.

  Therefore, we have proved there are projection $p\in B_{+}$ and self-adjoint element $y\in (1-p)A(1-p)$ with $\|y\|\leq1$ such that $\| y-x\|<\frac{\varepsilon}{3}$.

\vspace{10pt}
  Step $\uppercase\expandafter{\romannumeral2}$, we approximate $y$ by the sum of two elements: one is an invertible self-adjoint element in $(1-p)B(1-p)$, and the other one could be decomposed into the difference of two positive elements which are Cuntz subequivalent to a projection smaller than $p$.

   Since $B$ is unital simple infinite dimensional C*-algebra with property (SP), then there are nonzero mutually orthogonal projections $p_{1},\,p_{2},\,p_{3},\,p_{4},\,p_{5}$ such that $p_{1}\sim_{B} p_{2}\sim_{B} p_{3}\sim_{B} p_{4},\,p_{1}+p_{2}+p_{3}+p_{4}+p_{5}= p$.
   By lemma \ref{property}, there exists a nonzero projection $q_{1}\in (1-p)B(1-p)$ such that $q_{1}\precsim_{B}p_{1}$.

  Let $y_{+}=\frac{|y|+y}{2},\,y_{-}=\frac{|y|-y}{2}$ where $|y|=(y^{2})^{\frac{1}{2}}$, then $ y_{+},\,y_{-}\in (1-p)A(1-p)_{+}$ with $y=y_{+}-y_{-}$, $\| y_{+}\|\leq 1$ and $\| y_{-}\|\leq 1$.

   By lemma \ref{squares}, there is $\delta_{0}>0$ such that $\|e^{\frac{1}{2}}f-fe^{\frac{1}{2}}\|<\frac{\varepsilon}{24}$ whenever $e,\,f\in A$, $0\leq e\leq1,\, \| f\|\leq1,\, \| ef-fe\|<\delta_{0}$. By lemma \ref{dist}, there is $\delta_{1}>0$ such that if $w\in A_{+}$, $e\in B$ with $\|w\|\leq1,\,0\leq e\leq1$, there exists $f\in B_{+}$ such that $\| ewe-f\|<\delta_{1}$, then there exists $c\in B_{+}$ such that $\| e^{\frac{1}{2}}we^{\frac{1}{2}}-c\|<\frac{\varepsilon}{24}$.

  By lemma \ref{projection}, for $\{y_{+},\,y_{-}\}\subset (1-p)A(1-p)_{+}$, $\delta=\min\{\delta_{0},\,\delta_{1}\}$, $q_{1}\in ((1-p)B(1-p))_{+}\setminus\{0\}$, since  $B\subset A$ is a centrally large subalgebra, then there is $z_{1},\,z_{2}\in (1-p)B(1-p)_{+},\,g\in (1-p)B(1-p)$ such that:

 (1) $0\leq g\leq 1$;

 (2) $\|(1-p-g)y_{+}(1-p-g)-z_{1}\|<\delta,\,\|(1-p-g)y_{-}(1-p-g)-z_{2}\|<\delta$;

(3) $g\precsim_{(1-p)B(1-p)}q_{1}$;

(4)  $\|gy_{+}-y_{+}g\|<\delta,\,\| gy_{-}-y_{-}g\|<\delta$.

  By the choice of $\delta$ and (4), we have
  \begin{align}
  &\|g^{\frac{1}{2}}y_{+}-y_{+}g^{\frac{1}{2}}\|<\frac{\varepsilon}{24},\notag\\
   & \| g^{\frac{1}{2}}y_{-}-y_{-}g^{\frac{1}{2}}\|<\frac{\varepsilon}{24},\notag\\
   & \|(1-p-g)^{\frac{1}{2}}y_{+}-y_{+}(1-p-g)^{\frac{1}{2}}\|<\frac{\varepsilon}{24},\notag\\
   & \| (1-p-g)^{\frac{1}{2}}y_{-}-y_{-}(1-p-g)^{\frac{1}{2}}\|<\frac{\varepsilon}{24}\notag.
\end{align}

   Set $y_{1}=(1-p-g)^{\frac{1}{2}}y(1-p-g)^{\frac{1}{2}},\, y_{2}=g^{\frac{1}{2}}yg^{\frac{1}{2}}$, then we get
 \begin{align}
  &\| y-(y_{1}+y_{2})\|\notag\\
  &\leq\| y_{+}-(1-p-g)^{\frac{1}{2}}y_{+}(1-p-g)^{\frac{1}{2}}-g^{\frac{1}{2}}y_{+}g^{\frac{1}{2}}\|\notag\\
  &+\| y_{-}-(1-p-g)^{\frac{1}{2}}y_{-}(1-p-g)^{\frac{1}{2}}-g^{\frac{1}{2}}y_{-}g^{\frac{1}{2}}\| \notag\\
  &<\frac{\varepsilon}{24}+\frac{\varepsilon}{24}+\frac{\varepsilon}{24}+\frac{\varepsilon}{24}=\frac{\varepsilon}{6}\notag.
 \end{align}

  By the choice of $\delta$ and (2), then there are $w_{1}^{'},w_{2}^{'} \in ((1-p)B(1-p))_{+}$ such that
   $$\|(1-p-g)^{\frac{1}{2}}y_{+}(1-p-g)^{\frac{1}{2}}-w_{1}^{'}\|<\frac{\varepsilon}{24},$$
  $$ \|(1-p-g)^{\frac{1}{2}}y_{-}(1-p-g)^{\frac{1}{2}}-w_{2}^{'} \|<\frac{\varepsilon}{24}.$$

    Then $w_{1}^{'}-w_{2}^{'}\in (1-p)B(1-p)$ is a self-adjoint element. Since $RR((1-p)B(1-p))=0$, there is an invertible self-adjoint element $w_{1}\in (1-p)B(1-p)$ such that $\| w_{1}-(w_{1}^{'}-w_{2}^{'})\|<\frac{\varepsilon}{12}$. Then

     \begin{align}
     \| y_{1}-w_{1}\|&\leq\|y_{1}-(w_{1}^{'}-w_{2}^{'})\|+\|(w_{1}^{'}-w_{2}^{'})-w_{1}\|\notag\\
     &<\|(1-p-g)^{\frac{1}{2}}y_{+}(1-p-g)^{\frac{1}{2}}-(1-p-g)^{\frac{1}{2}}y_{-}(1-p-g)^{\frac{1}{2}}-(w_{1}^{'}-
     w_{1}^{'})\|+\frac{\varepsilon}{12}\notag\\
     &<\frac{\varepsilon}{24}+\frac{\varepsilon}{24}+\frac{\varepsilon}{12}=\frac{\varepsilon}{6}\notag.
  \end{align}

    Thus $\|(y_{1}+y_{2})-(w_{1}+y_{2})\|<\frac{\varepsilon}{6}$, so $\|y-(w_{1}+y_{2})\|<\frac{\varepsilon}{3}.$

And
    $$y_{2}=g^{\frac{1}{2}}yg^{\frac{1}{2}}=g^{\frac{1}{2}}y_{+}g^{\frac{1}{2}}-g^{\frac{1}{2}}y_{-}g^{\frac{1}{2}},$$
    $$g^{\frac{1}{2}}y_{+}g^{\frac{1}{2}}\leq \| y_{+}\|g\precsim_{B}q_{1}\precsim_{B} p_{1},$$
     $$g^{\frac{1}{2}}y_{-}g^{\frac{1}{2}}\leq \| y_{-}\|g\precsim_{B}q_{1}\precsim_{B} p_{1}.$$

\vspace{10pt}
   Step $\uppercase\expandafter{\romannumeral3}$, by matrix decomposition and some twirls, we get an invertible self-adjoint element in $A$ that could approximate to $w_{1}+y_{2}$.

    Since $g^{\frac{1}{2}}y_{+}g^{\frac{1}{2}}\precsim _{B}p_{1}$, $g^{\frac{1}{2}}y_{-}g^{\frac{1}{2}}\precsim_{B} p_{1}\sim p_{3}$, there exist $r_{j},c_{k}\in B$ such that
    $$\| r_{j}p_{1}r_{j}^{*}-g^{\frac{1}{2}}y_{+}g^{\frac{1}{2}}\|<\frac{\varepsilon}{12},\|\,\,\, c_{k}p_{3}c_{k}^{*}-g^{\frac{1}{2}}y_{-}g^{\frac{1}{2}}\|<\frac{\varepsilon}{12}.$$

  Then
  $$ \|r_{j}p_{1}r_{j}^{*}-c_{k}p_{3}c_{k}^{*}-y_{2}\|<\frac{\varepsilon}{6},$$
  $$\| p_{1}r_{j}^{*}\|=\|r_{j}p_{1}\|\leq\sqrt{1+\frac{\varepsilon}{12}}<2,$$
  $$\| p_{3}c_{k}^{*}\|=\|c_{k}p_{3}\|\leq\sqrt{1+\frac{\varepsilon}{12}}<2.$$

    Since $p_{1}\sim_{B} p_{2},\,p_{3}\sim_{B} p_{4}$ , there are $v_{1},\,v_{2} \in B$ such that
    $$p_{1}=v_{1}^{*}v_{1},\, p_{2}=v_{1}v_{1}^{*},\,p_{3}=v_{2}^{*}v_{2},\, p_{4}=v_{2}v_{2}^{*}.$$

   Let
   $$r=\left(
         \begin{array}{cccccc}
           1 & (1-p)r_{j}p_{1} & 0 & (1-p)c_{k}p_{3} & 0 & 0 \\
          0 & \frac{\varepsilon}{6} & 0 & 0 & 0 & 0 \\
           0 & 0 & \frac{\varepsilon}{6}& 0 & 0 & 0 \\
           0 & 0 & 0 &\frac{\varepsilon}{6} & 0 & 0 \\
           0 & 0 & 0 & 0 & \frac{\varepsilon}{6} & 0 \\
           0 & 0 & 0 & 0 & 0 & \frac{\varepsilon}{6} \\
         \end{array}
       \right),$$
and
   $$w=\left(
         \begin{array}{cccccc}
           w_{1} & 0 & 0 & 0 & 0 & 0 \\
           0 & p_{1} & \frac{\varepsilon}{6}v_{1}^{*} & 0 & 0 & 0 \\
           0 & \frac{\varepsilon}{6}v_{1} & 0 & 0 & 0 & 0 \\
           0 & 0 & 0 & -p_{3} & \frac{\varepsilon}{6}v_{2}^{*} & 0 \\
           0 & 0 & 0 & \frac{\varepsilon}{6}v_{2} & 0 & 0 \\
           0& 0 & 0 & 0 & 0 & \frac{\varepsilon}{6} \\
         \end{array}
       \right).$$

       We obtain
       $$w^{-1}=\left(
         \begin{array}{cccccc}
           w_{1}^{-1} & 0 & 0 & 0 & 0 & 0 \\
           0 & 0 & \frac{6}{\varepsilon}v_{1}^{*} & 0 & 0 & 0 \\
           0 & \frac{6}{\varepsilon}v_{1} &-\frac{36}{\varepsilon^{2}}v_{1}p_{1}v_{1}^{*} & 0 & 0 & 0 \\
           0 & 0 & 0 & 0 & \frac{6}{\varepsilon}v_{2}^{*} & 0 \\
           0 & 0 & 0 & \frac{6}{\varepsilon}v_{2} & \frac{36}{\varepsilon^{2}}v_{2}p_{3}v_{2}^{*} & 0 \\
           0& 0 & 0 & 0 & 0 & \frac{6}{\varepsilon} \\
         \end{array}
       \right).$$

       Then $r$ is invertible in $A$ and $w$ is invertible and self-adjoint in $A$. Moreover,

  $rwr^{*}=$

      $$\left(\begin{smallmatrix} w_{1}+(1-p)(r_{j}p_{1}r_{j}^{*}-c_{k}p_{3}c_{k}^{*})(1-p)  & \frac{\varepsilon}{6}(1-p)r_{j}p_{1}  & (\frac{\varepsilon}{6})^{2}(1-p)r_{j}p_{1}v_{1}^{*}  & -\frac{\varepsilon}{6}(1-p)c_{k}p_{3} & (\frac{\varepsilon}{6})^{2}(1-p)c_{k}p_{3}v_{2}^{*}  & 0  \\ \frac{\varepsilon}{6} p_{1}r_{j}^{*}(1-p) & (\frac{\varepsilon}{6})^{2}\!p_{1}  & (\frac{\varepsilon}{6})^{3} \!v_{1}^{*} & 0 & 0  & 0 \\(\frac{\varepsilon}{6})^{2} \!v_{1}p_{1}r_{j}^{*} (1-p) & (\frac{\varepsilon}{6})^{3} v_{1} & 0 & 0  & 0  & 0 \\
                     -\frac{\varepsilon}{6}p_{3}c_{k}^{*}(1-p)  & 0  & 0  & -(\frac{\varepsilon}{6})^{2}p_{3}  & (\frac{\varepsilon}{6})^{3}v_{2}^{*} & 0  \\
                     (\frac{\varepsilon}{6})^{2}\!v_{2}p_{3}c_{k}^{*}(1-p) & 0 & 0  & (\frac{\varepsilon}{6})^{3}\!v_{2}  & 0 & 0  \\
                     0 & 0 & 0  & 0  & 0 & (\frac{\varepsilon}{6})^{3}  \end{smallmatrix} \right),$$

       $\| rwr^{*}-(w_{1}+y_{2})\|$
      $$ =\left\Vert\left(
                   \begin{smallmatrix}
                 (1-p)[r_{j}p_{1}r_{j}^{*}-c_{k}p_{3}c_{k}^{*}-\!y_{2}](1-p)  & \frac{\varepsilon}{6}(1-p)r_{j}p_{1}  & (\frac{\varepsilon}{6})^{2}(1-p)r_{j}p_{1}v_{1}^{*}  & -\frac{\varepsilon}{6} (1-p)c_{k}p_{3} & (\frac{\varepsilon}{6})^{2}(1-p)c_{k}p_{3}v_{2}^{*}  & 0  \\
                     \frac{\varepsilon}{6} p_{1}r_{j}^{*}(1-p) & (\frac{\varepsilon}{6})^{2}p_{1}  & (\frac{\varepsilon}{6})^{3} v_{1}^{*} & 0 & 0  & 0 \\
                     (\frac{\varepsilon}{6})^{2} \!v_{1}p_{1}r_{j}^{*} (1-p) & (\frac{\varepsilon}{6})^{3} \!v_{1} & 0 & 0  & 0  & 0 \\
                     -\frac{\varepsilon}{6}\!p_{3}c_{k}^{*}(1-p)  & 0  & 0  & -(\frac{\varepsilon}{6})^{2}\!p_{3}  & (\frac{\varepsilon}{6})^{3}v_{2}^{*} & 0  \\
                     (\frac{\varepsilon}{6})^{2}\!v_{2}p_{3}c_{k}^{*}(1-p) & 0 & 0  & (\frac{\varepsilon}{6})^{3}\!v_{2}  & 0 & 0  \\
                     0 & 0 & 0  & 0  & 0 & (\frac{\varepsilon}{6})^{3}
                   \end{smallmatrix}
                 \right)\right\Vert$$

       $<\frac{\varepsilon}{3}.$

       Hence $rwr^{*}$ is an invertible self-adjoint element in $A$, and
 \begin{align}
  &\| x-rwr^{*}\|\notag\\
   &\leq \| x-y\|+\|y-(w_{1}+y_{2})\|+\|(w_{1}+y_{2})-rwr^{*}\|\notag\\
   & <\frac{\varepsilon}{3}+\frac{\varepsilon}{3}+\frac{\varepsilon}{3}=\varepsilon\notag.
\end{align}

\end{proof}

\section{Weak Comparison}

 Local weak comparison and weak comparison are first introduced by Kirchberg and R{\o}rdam in \cite{E. K1}. In this section, we give the permanence of local weak comparison and weak comparison.

\begin{definition}\label{weak comparison}
Let $A$ be a unital, simple, stably finite C*-algebra. We say $A$ has $local\,\, weak\,\,comparison$, if there is a constant $\gamma(A)\in[1,\infty)$ such that the following holds. For all positive elements $a$ and $b$ in $A$: If
\[\gamma(A)\cdot \sup_{\tau\in QT(A)} d_{\tau}(a)<\inf_{\tau\in QT(A)} d_{\tau}(b),\]
then $\langle a\rangle\leq\langle b\rangle$ in the Cuntz semigroup Cu$(A)$ of $A$.

  If $M_{n}(A)$ has local weak comparison for all $n$, and $\sup_{n}\gamma(M_{n}(A))<\infty$, then we say
that $A$ has $weak\,\, comparison$.
\end{definition}

First we show that if $A $ has  local weak comparison, then its large subalgebra $B$ has local weak comparison.
Since we need to discuss $\gamma(A)\cdot \sup_{\tau\in QT(A)} d_{\tau}(a)<\inf_{\tau\in QT(A)} d_{\tau}(b)$, so the following form of Dini's theorem is useful (see the proof in lemma 6.13 of \cite{N. C. P2}).

\begin{lem}\label{continuous 1}\cite[Lemma 6.13]{N. C. P2}
 Let $X$ be a compact Hausdorff space. Let $(f_{n})_{n\in \mathbb{Z}_{+}}$ be a sequence of lower semicontinuous functions $f_{n}:X\rightarrow \mathbb{R}\bigcup\{\infty\}$ such that $f_{1}(x)\leq f_{2}(x)\leq...,$ for all $x\in X$ and let $g:X\rightarrow\mathbb{R}$ be a continuous function such that $g(x)<lim_{n\rightarrow\infty}f_{n}(x)$ for all $x\in X$. Then there is $n\in \mathbb{Z}_{>0}$ such that $f_{n}(x)>g(x)$ for all $x\in X$.
 \end{lem}

 \begin{lem}\label{smaller}

 Let $A$ be a unital, simple, stably finite C*-algebra. For all positive elements $a$ and $b$ in $A$, $\gamma \in \mathbb{R}$ such that
\[\gamma\cdot \sup_{\tau\in QT(A)} d_{\tau}(a)<\inf_{\tau\in QT(A)} d_{\tau}(b).\]

 Then there exists $\varepsilon>0$ such that
$$\gamma\cdot \sup_{\tau\in QT(A)} d_{\tau}(a)<\inf_{\tau\in QT(A)} d_{\tau}((b-\varepsilon)_{+}).$$
 \end{lem}

 \begin{proof}
  Since $\gamma\cdot \sup_{\tau\in QT(A)} d_{\tau}(a)<\inf_{\tau\in QT(A)} d_{\tau}(b)$, there exists $\delta>0$ such that
  \[\gamma\cdot \sup_{\tau\in QT(A)} d_{\tau}(a)+\delta<\inf_{\tau\in QT(A)} d_{\tau}(b).\]

 Define $f_{n}:QT(A)\rightarrow[0,\infty]$ by $f_{n}(\tau)=d_{\tau}(b-\frac{1}{n})_{+}$ for $\tau \in QT(A)$ and $n\in \mathbb{Z}_{>0}$, and
define $g:QT(A)\rightarrow[0,\infty)$ by $g(\tau)=\gamma\cdot \sup_{\tau\in QT(A)} d_{\tau}(a)+\delta$ for $\tau \in QT(A)$.
Then $\{f_{n}\}_{n\in \mathbb{Z}_{>0}}$ are lower semicontinuous functions by lemma \ref{lower}, $g$ is a continuous function, and

 $$g(\tau)<\inf_{\tau\in QT(A)} d_{\tau}(b)\leq d_{\tau}(b)=\lim_{n\to\infty}f_{n}(\tau).$$

 So by lemma \ref{continuous 1}, we can get a $n\in \mathbb{Z}_{>0}$ such that $f_{n}(\tau)>g(\tau)$ for all $\tau \in QT(A)$. That is,
 $$\gamma\cdot \sup_{\tau\in QT(A)} d_{\tau}(a)+\delta<d_{\tau}((b-\frac{1}{n})_{+}).$$

 Let $\varepsilon=\frac{1}{n}$. It follows that
 \[\gamma\cdot \sup_{\tau\in QT(A)} d_{\tau}(a)<\gamma\cdot \sup_{\tau\in QT(A)} d_{\tau}(a)+\delta\leq\inf_{\tau\in QT(A)} d_{\tau}((b-\varepsilon)_{+}).\]

 \end{proof}
\begin{lem}\label{type 1}\cite[Lemma 2.7]{N. C. P2}
Let $A$ be a simple infinite dimensional C*-algebra which is not of type $\uppercase\expandafter{\romannumeral1}$. Let $b\in A_{+}\setminus\{0\}$, $\varepsilon>0$, and $n\in \mathbb{Z}_{+}$. Then there are $c\in A_{+}, \,y\in A_{+}\setminus\{0\}$ such that
$$n\langle(b-\varepsilon)_{+}\rangle\leq(n+1)\langle c\rangle, \,\,\,\langle c\rangle+\langle y\rangle\leq\langle b\rangle$$
in $W(A)$.
\end{lem}

\begin{thrm}\label{A implies B}
Let $A$ be an infinite dimensional stably finite simple separable unital C*-algebra, and let $B\subset A$ be a large subalgebra.  If $A $ has  local weak comparison, then $B$ has local weak comparison.
\end{thrm}

\begin{proof}

Let $\gamma(B)=\gamma(A)$,  for all positive elements $a$ and $b$ in $B$: If

\[\gamma(A)\cdot \sup_{\tau\in QT(B)} d_{\tau}(a)<\inf_{\tau\in QT(B)} d_{\tau}(b),\]
then we show will $\langle a\rangle\leq\langle b\rangle$ in the Cuntz semigroup Cu$(B)$ of $B$. Thus we will prove that $(a-\varepsilon)_{+}\precsim_{B}b$ for all $\varepsilon>0$ by lemma \ref{fundamental}(4).

By lemma \ref{smaller}, there exists $\varepsilon_{0}>0$ such that
$$r=\inf_{\tau\in QT(B)} d_{\tau}((b-\varepsilon_{0})_{+})-\gamma(A)\cdot \sup_{\tau\in QT(B)} d_{\tau}(a)>0.$$

Let $n>\frac{\inf_{\tau\in QT(B)} d_{\tau}((b-\varepsilon_{0})_{+})}{r}-1$ and $\varepsilon_{0}>0$. According to proposition 5.2 and 5.5 in \cite{N. C. P2}, $B$ is simple infinite dimension C*-algebra, so that it is not type  $\uppercase\expandafter{\romannumeral1}$. By lemma \ref{type 1}, there exist  $c\in B_{+},\, y\in B_{+}\setminus\{0\}$ such that
\begin{align}
n\langle(b-\varepsilon_{0})_{+}\rangle\leq(n+1)\langle c\rangle, \,\,\,\langle c\rangle+\langle y\rangle\leq\langle b\rangle
\end{align}
in $W(B)$.

Then $\frac{n}{n+1} d_{\tau}((b-\varepsilon_{0})_{+})\leq d_{\tau}(c)$ for all $\tau\in QT(B)$ by notation \ref{quasitraces}. Hence

$$ \inf_{\tau\in QT(B)} d_{\tau}(c)\geq \frac{n}{n+1} \inf_{\tau\in QT(B)} d_{\tau}((b-\varepsilon_{0})_{+}),$$

$$\inf_{\tau\in QT(B)} d_{\tau}(c)-\gamma(A)\cdot \sup_{\tau\in QT(B)} d_{\tau}(a)\geq \frac{n}{n+1} \inf_{\tau\in QT(B)} d_{\tau}((b-\varepsilon_{0})_{+})-\gamma(A)\cdot \sup_{\tau\in QT(B)} d_{\tau}(a).$$

By the choice of $n$, we have
\begin{align}
  &\frac{n}{n+1} \inf_{\tau\in QT(B)} d_{\tau}((b-\varepsilon_{0})_{+})-\gamma(A)\cdot \sup_{\tau\in QT(B)} d_{\tau}(a)\notag  \\
  &=\inf_{\tau\in QT(B)} d_{\tau}((b-\varepsilon_{0})_{+})-\gamma(A)\cdot \sup_{\tau\in QT(B)} d_{\tau}(a)-\frac{1}{n+1} \inf_{\tau\in QT(B)} d_{\tau}((b-\varepsilon_{0})_{+}) \notag\\
  &=r-\frac{1}{n+1} \inf_{\tau\in QT(B)} d_{\tau}((b-\varepsilon_{0})_{+})>0\notag.
\end{align}

Then it follows that
\[\gamma(A)\cdot \sup_{\tau\in QT(B)} d_{\tau}(a)<\inf_{\tau\in QT(B)} d_{\tau}(c).\]

Since $A$ has local weak comparison, we have $a\precsim _{A}c$. So there is $v\in A$ such that $\|vcv^{*}-a\|<\varepsilon$.

Since $B\subseteq A$ is a large subalgebra, then for $F=\{v\}$, $\frac{\varepsilon}{4\|c\|\|v\|+1}$ and $y\in B$, there are $v_{0}\in A_{+}$ and $g\in B$ such that:

(1) $0\leq g\leq1;$

(2) $\parallel v-v_{0}\parallel<\frac{\varepsilon}{4\|c\|\|v\|+1};$

(3) $\parallel v_{0}\parallel\leq\parallel v\parallel;$

(4) $(1-g)v_{0}\in B;$

(5) $g\precsim_{B}y.$

Then we get $\|v_{0}cv_{0}-vcv\|<\frac{\varepsilon}{2}$, and so
$$\|(1-g)v_{0}c[(1-g)v_{0}]^{*}-(1-g)a(1-g)\|<\varepsilon.$$

By lemma \ref{fundamental}(3), we have
\begin{align}
 [(1-g)a(1-g)-\varepsilon]_{+}\precsim_{B}(1-g)v_{0}c[(1-g)v_{0}]^{*}\precsim_{B}c.
\end{align}

So $$(a-\varepsilon)_{+}\precsim_{B}[(1-g)a(1-g)-\varepsilon]_{+}\oplus g\precsim_{B}c\oplus g\precsim_{B}c\oplus y\precsim_{B}b.$$

Use lemma \ref{fundamental}(6) at the first step,  (4.2) at the second step, (5) at the third step, (4.1) at the forth step. Therefore, we prove $\langle a\rangle\leq\langle b\rangle$ in the Cuntz semigroup Cu$(B)$ of $B$.

\end{proof}

Next we prove the converse direction. However, we can not get $A$ has  local weak comparison if its large subalgebra $B$ has local weak comparison, we obtain that $A$ has  local weak comparison  when $M_{2}(B)$ has local weak comparison.

The following lemma could be found in \cite{N. C. P2} and it said that we could choose an element very small in the sense of quasitrace in a simple unital infinite dimensional $C$*-algebra.

\begin{lem}\label{small}\cite[Corollary 2.5]{N. C. P2}
Let $A$ be a simple unital infinite dimensional C*-algebra. Then for every $\varepsilon>0$, there is $a\in A_{+}\setminus\{0\}$ such that $d_{\tau}(a)<\varepsilon$ for all $\tau\in QT(A)$.
\end{lem}

\begin{thrm}\label{B implies A}
 Let $A$ be an infinite dimensional stably finite simple separable unital C*-algebra, and let $B\subset A$ be a large subalgebra.  If $M_{2}(B)$ has local weak comparison, then $A$ has  local weak comparison.
\end{thrm}

\begin{proof}

Let $\gamma(A)=\gamma(M_{2}(B))$,  for all positive elements $a$ and $b$ in $A$: If
\[\gamma(M_{2}(B))\cdot \sup_{\tau\in QT(A)} d_{\tau}(a)<\inf_{\tau\in QT(A)} d_{\tau}(b),\]
then we will show $\langle a\rangle\leq\langle b\rangle$ in the Cuntz semigroup Cu$(A)$ of $A$. Next, we will prove $(a-\varepsilon)_{+}\precsim_{A}b$ for any $\varepsilon>0$.

By lemma \ref{smaller}, there exists $\varepsilon_{0}$ such that

$$\gamma(M_{2}(B))\cdot \sup_{\tau\in QT(A)} d_{\tau}(a)<\inf_{\tau\in QT(A)} d_{\tau}((b-\varepsilon_{0})_{+}).$$

Let $r=\inf_{\tau\in QT(A)} d_{\tau}((b-\varepsilon_{0})_{+})-\gamma(M_{2}(B))\cdot \sup_{\tau\in QT(A)} d_{\tau}(a)$.

Since $B\subseteq A$ is a large subalgebra, then for $F=\{a,b\}$, and  $0<\varepsilon_{1}<\min\{\frac{\varepsilon}{2},\frac{\varepsilon_{0}}{2},\frac{r}{2}\}$, by lemma \ref{small}, there is $x\in A_{+}\setminus\{0\}$ such that $d_{\tau}(x)<\frac{\varepsilon_{1}}{\gamma(M_{2}(B))}$.
Then there are $a_{0},\,b_{0}\in A_{+}$ and $g\in B$ such that:

(1) $0\leq g\leq1;$

(2) $\parallel a-a_{0}\parallel<\varepsilon_{1},\,\parallel b-b_{0}\parallel<\varepsilon_{1};$

(3) $(1-g)a_{0}(1-g)\in B,\,(1-g)b_{0}(1-g)\in B;$

(4) $g\precsim_{A}x.$

By the choice of $\varepsilon_{1}$, we have
\begin{align}
&\gamma(M_{2}(B))d_{\tau}(g)<\varepsilon_{1},\\
&\|(1-g)a_{0}(1-g)-(1-g)a(1-g)\|<\frac{\varepsilon}{2},\\
&\inf_{\tau\in QT(A)} d_{\tau}((b-\varepsilon_{0})_{+})-\gamma(M_{2}(B))\cdot \sup_{\tau\in QT(A)} d_{\tau}(a)-2\varepsilon_{1}>0,\\
&\|(1-g)b(1-g)-(1-g)b_{0}(1-g)\|<\varepsilon_{0}-\varepsilon_{1}.
\end{align}

Then by (4.3) at the second step. (4.4) and lemma \ref{fundamental}(3) at the third step, we can obtain
\begin{align}
&\gamma(M_{2}(B))\cdot \sup_{\tau\in QT(A)} d_{\tau}(((1-g)a_{0}(1-g)-\frac{\varepsilon}{2})_{+}\oplus g)\notag\\
&\leq\gamma(M_{2}(B))\cdot \sup_{\tau\in QT(A)} d_{\tau}(((1-g)a_{0}(1-g)-\frac{\varepsilon}{2})_{+})+\gamma(M_{2}(B))\cdot \sup_{\tau\in QT(A)} d_{\tau}(g)\notag\\
&\leq\gamma(M_{2}(B))\cdot \sup_{\tau\in QT(A)} d_{\tau}(((1-g)a_{0}(1-g)-\frac{\varepsilon}{2})_{+})+\varepsilon_{1}\notag\\
&\leq\gamma(M_{2}(B))\cdot \sup_{\tau\in QT(A)} d_{\tau}(((1-g)a(1-g))+\varepsilon_{1}\notag\\
&\leq\gamma(M_{2}(B))\cdot \sup_{\tau\in QT(A)} d_{\tau}(a)+\varepsilon_{1}\notag.
   \end{align}
By (4.5) at the first step, lemma \ref{fundamental}(5) with $\lambda=\varepsilon_{0}-\varepsilon_{1}$ at the second step, lemma \ref{fundamental}(6) at the third step, $d_{\tau}(g)\leq d_{\tau}(x)<\frac{\varepsilon_{1}}{\gamma(M_{2}(B))}\leq \varepsilon_{1}$ at the forth step, we get
\begin{align}
&\gamma(M_{2}(B))\cdot \sup_{\tau\in QT(A)} d_{\tau}(a)+\varepsilon_{1}+\varepsilon_{1}\notag\\
&<\inf_{\tau\in QT(A)} d_{\tau}((b-\varepsilon_{0})_{+}), \notag\\
&\leq \inf_{\tau\in QT(A)} d_{\tau}((b_{0}-(\varepsilon_{0}-\varepsilon_{1}))_{+})\notag\\
&\leq \inf_{\tau\in QT(A)} d_{\tau}(((1-g)b_{0}(1-g)-(\varepsilon_{0}-\varepsilon_{1}))_{+})+\inf_{\tau\in QT(A)} d_{\tau}(g)\notag\\
&\leq \inf_{\tau\in QT(A)} d_{\tau}(((1-g)b_{0}(1-g)-(\varepsilon_{0}-\varepsilon_{1}))_{+})+\varepsilon_{1}\notag.
\end{align}
Thus
$$\gamma(M_{2}(B))\cdot \sup_{\tau\in QT(A)} d_{\tau}(a)+\varepsilon_{1}< \inf_{\tau\in QT(A)} d_{\tau}(((1-g)b_{0}(1-g)-(\varepsilon_{0}-\varepsilon_{1}))_{+}).$$
It follows that
$$ \gamma(M_{2}(B))\cdot \sup_{\tau\in QT(A)} d_{\tau}(((1-g)a_{0}(1-g)-\frac{\varepsilon}{2})_{+}\oplus g)< \inf_{\tau\in QT(A)} d_{\tau}(((1-g)b_{0}(1-g)-(\varepsilon_{0}-\varepsilon_{1}))_{+}).$$

Since $((1-g)a_{0}(1-g)-\frac{\varepsilon}{2})_{+}\oplus g,((1-g)b_{0}(1-g)-(\varepsilon_{0}-\varepsilon_{1}))_{+}\in M_{2}(B)$, and $M_{2}(B)$ has local weak comparison, then
\begin{align}
((1-g)a_{0}(1-g)-\frac{\varepsilon}{2})_{+}\oplus g\precsim_{B}((1-g)b_{0}(1-g)-(\varepsilon_{0}-\varepsilon_{1}))_{+}.
\end{align}
Thus
\begin{align}
  (a-\varepsilon)_{+}& \precsim_{A}(a_{0}-\frac{\varepsilon}{2})_{+}\notag \\
   & \precsim_{A}((1-g)a_{0}(1-g)-\frac{\varepsilon}{2})_{+}\oplus g\notag\\
   & \precsim_{B}[(1-g)b_{0}(1-g)-(\varepsilon_{0}-\varepsilon_{1})]_{+}\notag\\
   & \precsim_{A}(1-g)b(1-g)\notag\\
   & \precsim_{A} b\notag.
\end{align}
Use lemma \ref{fundamental}(5) with $\lambda=\frac{\varepsilon}{2}$ at the first step, lemma \ref{fundamental}(6) at the second step, (4.7) at the third step, (4.6) with lemma \ref{fundamental}(3) at the forth step.

Therefore, we prove $\langle a\rangle\leq\langle b\rangle$ in the Cuntz semigroup Cu$(A)$ of $A$.
\end{proof}

With the theorems before, we consider the permanence of weak comparison.

\begin{cor} Let $A$ be an infinite dimensional stably finite simple separable unital C*-algebra, and let $B\subset A$ be a large subalgebra. Then $A$ has  weak comparison if and only if $B$ has weak comparison.
\end{cor}

\begin{proof}
First we suppose that $B$ has weak comparison, then $M_{n}(B)$ has local weak comparison for all $n$, and $\sup_{n}\gamma(M_{n}(B))<\infty$.
By theorem \ref{B implies A}, since $M_{2n}(B)$ has local weak comparison, we have $M_{n}(A)$ has local weak comparison, and $\gamma(M_{n}(A))=\gamma(M_{2n}(B))$. So $\sup_{n}\gamma(M_{n}(A))<\infty$. It follows that $A$ has weak comparison.

For the other direction,  since $A$ has weak comparison, then $M_{n}(A)$ has local weak comparison for all $n$, and $\sup_{n}\gamma(M_{n}(A))<\infty$. By theorem \ref{A implies B}, then $M_{n}(B)$ has local weak comparison for all $n$, and $\gamma(M_{n}(B))=\gamma(M_{n}(A))$, so $\sup_{n}\gamma(M_{n}(B))<\infty$. That is, $B$ has weak comparison.

\end{proof}
\noindent {\bfseries Acknowledgement}

~\\This work was partially supported by National Natural Science Foundation of China [Grant No. 11871375].


\begin{thebibliography}{}
\bibitem{I. F. P1} I. F. Putnam, \emph{The C*-algebras associated with mimimal homeomorphism of the Cantor set,}  Pac. J. Math. 136.2(1989): 329-353.
\bibitem{Q. L1} Q. Lin and N. C. Phillips, \emph{ordered K-theory for C*-algebras of minimal homeomorphisms,} Operator Algebras and Operator Theory, L. Ge, etc., eds., Contemp. Math. 228(1998): 289-314.
\bibitem{G. A. E1} G. A. Elliott and Z. Niu, \emph{The C*-algebras of minimal homeomorphism of zero mean dimension,} preprint (arXiv:1406.2382v2 [math.OA]).
\bibitem{N. C. P1} N. C. Phillips, \emph{Cancellation and stable rank for direct limits of recursive subhomogeneous algebras,} T. Am. Math. Soc. 359.10(2007): 4625-4652.

\bibitem{N. C. P2} N. C. Phillips, \emph{Large subalgebras,} preprint(arXiv: 1408.5546v1[math. OA]).

\bibitem{N. C. P3} D. Archey and N. C. Phillips, \emph{Permanence of stable rank one for centrally large subalgebras and crossed products by minimal homeomorphisms,} preprint(arXiv: 1505.00725v2v1[math. OA]).

\bibitem{N. C. P4} D. Archey, J. Buck and N. C. Phillips, \emph{Centrally large subalgebras and tracial $\mathcal{Z}$-stability,} Int. Math. Res. Notices 2018.6(2017): 1857-1877.

\bibitem{N. C. P5} N. C. Phillips, \emph{Large subagebras and the structure of crossed products(draft),} (2015).

\bibitem{A. T1} W. Winter, \emph{Decomposition rank and $\mathcal{Z}$-stability,} Invent. Math. 179.2(2010): 229-301.

\bibitem{E. O1} E. Ortega, F. Perera, and M. R{\o}rdam, \emph{The Corona Factorization property, Stability, and the
Cuntz semigroup of a C*-algebra,} Int. Math. Res. Notices 2012.1(2012): 34-66.

 \bibitem{W. W1} W. Winter, \emph{Nuclear dimension and $\mathcal{Z}$-stability of pure C*-algebras,} Invent. Math. 187.2(2012): 259-342.

 \bibitem{E. K1} E. Kirchberg and M. R{\o}rdam, \emph{Central sequence $C$*-algebras and tensorial absorption of the Jiang-Su algebra,} J. Reine Angew. Math. 695(2014): 175-214

 \bibitem{Q. F1} Q. Z. Fan, X. C. Fang and X. Zhao, \emph{the comparison properties of large subalgebra are inherited,} prepared.

 \bibitem{B. B1} B. Blackadar and D. Handelman, \emph{Dimension functions and traces on C*-algebras,} J. Funct. Anal. 45.45(1982): 297-340.

 \bibitem{P. A1} P. Ara, F. Perera, and A. S. Toms, \emph{K-Theory for operator algebras. Classification of C*-
algebras,} In: Aspects of Operator Algebras and Applications. Contemp. Math., Providence, RI: Amer. Math. Soc. 534(2011): 1-71.

 \bibitem{G. A. E2} G. A. Elliott, L. Robert, and L. Santiago, \emph{The cone of lower semicontinuous traces on a
C*-algebra,} Am. J. Math. 133.4(2011): 969-1005.

\bibitem{H. X. L1}H. X. Lin, \emph{An introduction to the classification of amenable C*-algebras,} World Scientific, 2001.

\end{thebibliography}
\end{document}